\documentclass[12pt,a4paper]{amsart}
\usepackage{url,amscd,amssymb}

\newcommand{\A}{\mathcal{A}}
\newcommand{\D}{\mathcal{D}}
\newcommand{\DA}{\ensuremath{\D A}}
\newcommand{\DB}{\ensuremath{\D B}}
\newcommand{\Q}{\ensuremath{\mathbb{Q}}}
\newcommand{\Z}{\ensuremath{\mathbb{Z}}}
\newcommand{\F}{\ensuremath{\mathbb{F}}}
\newcommand{\R}{\ensuremath{\mathbb{R}}}
\newcommand{\Zp}{\ensuremath{\Z_{p}}}
\newcommand{\Zpl}{\ensuremath{\Z_{(p)}}}
\newcommand{\Zplu}{\ensuremath{\Zpl^{\times}}}
\newcommand{\Fp}{\ensuremath{\F_{p}}}
\newcommand{\Fpu}{\ensuremath{\Fp^{\times}}}
\newcommand{\padic}{$p$\nobreakdash-adic}
\newcommand{\padically}{\padic ally}
\newcommand{\plocal}{$p$\nobreakdash-local}
\newcommand{\Ktheory}{$K$\nobreakdash-theory}
\newcommand{\floor}[1]{\left\lfloor#1\right\rfloor}
\newcommand{\gp}[2]{\genfrac{[}{]}{0pt}{}{#1}{#2}}

\newcommand{\KzKp}{\ensuremath{K_{0}(K)_{(p)}}}
\newcommand{\tensor}{\otimes}
\newcommand{\MQ}{\ensuremath{M\tensor\Q}}
\newcommand{\AQ}{\ensuremath{A\tensor\Q}}
\newcommand{\set}[1]{\ensuremath{\{\,#1\,\}}}

\renewcommand{\Im}{\operatorname{Im}}
\newcommand{\Ker}{\operatorname{Ker}}
\newcommand{\Hom}{\operatorname{Hom}}
\newcommand{\Homz}{\Hom_{\Zpl}\!}
\newcommand{\Homc}{\Hom^{\operatorname{cts}}_{\Zpl}\!}
\newcommand{\Ext}{\operatorname{Ext}}
\newcommand{\Psiq}{\ensuremath{\Psi^{q}}}
\newcommand{\fbar}{\bar{f}}
\newcommand{\gbar}{\bar{g}}
\newcommand{\ft}{\tilde{f}}
\newcommand{\qh}{\hat{q}}
\newcommand{\Phih}{\hat{\Phi}}
\renewcommand{\_}{-}
\let\le\leqslant
\let\ge\geqslant
\let\phi\varphi
\let\epsilon\varepsilon

\theoremstyle{plain}
\newtheorem{thm}[equation]{Theorem}
\newtheorem{prop}[equation]{Proposition}
\newtheorem{cor}[equation]{Corollary}
\newtheorem{lemma}[equation]{Lemma}

\theoremstyle{definition}
\newtheorem{defn}[equation]{Definition}
\newtheorem{remark}[equation]{Remark}
\newtheorem{ex}[equation]{Example}

\numberwithin{equation}{section}

\begin{document}

\title[The discrete module category for $K$-theory operations]{The
discrete module category for the ring of $K$-theory operations}

\author[F. Clarke]{Francis Clarke}
\address{Department of Mathematics, University of Wales Swansea,
Swansea SA2 8PP, Wales}
\email{F.Clarke@Swansea.ac.uk}

\author[M. D. Crossley]{Martin Crossley}
\address{Department of Mathematics, University of Wales Swansea,
Swansea SA2 8PP, Wales}
\email{M.D.Crossley@Swansea.ac.uk}

\author[S. Whitehouse]{Sarah Whitehouse}
\address{Department of Pure Mathematics,
University of Sheffield,
Sheffield S3 7RH, England}
\email{S.Whitehouse@sheffield.ac.uk}

\begin{abstract}
    We study the category of discrete modules over the ring of degree
    zero stable operations in $p$-local complex $K$-theory.  We show
    that the $K_{(p)}$-homology of any space or spectrum is such a
    module, and that this category is isomorphic to a category defined
    by Bousfield and used in his work on the $K_{(p)}$-local stable
    homotopy category~\cite{Bousfield85}.  We also provide an
    alternative characterisation of discrete modules as locally
    finitely generated modules.
\end{abstract}



\date{$27^{\text{th}}$ February 2006}

\maketitle


\section{Introduction}
\label{intro}

An explicit description of the topological ring $A$ of degree zero
stable operations in \plocal\ complex \Ktheory\ was given
in~\cite{CCW05}.  Here we consider the category of `discrete modules'
over this ring.  We focus attention on \emph{discrete} modules because
the $K_{(p)}$-homology of any space or spectrum is such a module; see
Proposition~\ref{spectraarecontinuous}.

In Section~\ref{DiscreteModules} we recall some results
from~\cite{CCW05} about the ring $A$, define discrete $A$-modules, and
show how the category \DA\ of such modules is easily seen to be a
cocomplete abelian category.

Our first main result comes in Section~\ref{localfg}, where we provide
an interesting alternative characterisation of discrete $A$-modules:
an $A$-module is discrete if and only if it is locally finitely
generated.  We also show that the category \DA\ is isomorphic to the
category of comodules over the coalgebra \KzKp\ of which $A$ is the
dual.

The category \DA\ arose in disguised form in~\cite{Bousfield85}.
There Bousfield introduced a certain category $\A(p)$ as the first
step in his investigation of the $K_{(p)}$-local stable homotopy
category.  In Section~\ref{bousfieldscategory} we recall Bousfield's
(rather elaborate) definition, and we prove that his category is
isomorphic to the category of discrete $A$-modules.

This allows us to simplify and clarify in Section~\ref{cofreeobjs}
some constructions in Bousfield's work.  In particular, there is a
right adjoint to the forgetful functor from $\A(p)$ to the category of
\Zpl-modules.  For Bousfield, this functor has to be constructed in an
ad hoc fashion, separating cases.  In our context, it is revealed as
simply a continuous $\Hom$ functor.  We give a construction in this
language of a four-term exact sequence involving the right adjoint.

Bousfield's aim was to give an algebraic description of the
$K_{(p)}$-local stable homotopy category.  He succeeded at the level
of objects, and in Section~\ref{Bresult} we translate his main result
into our language of discrete $A$-modules.

Of course, \plocal\ \Ktheory\ splits as a sum of copies of the Adams
summand.  We have chosen to write the main body of this paper in the
non-split context, but very similar results hold in the split setting.
We record these in an appendix.

We note that a full algebraic description of the $K_{(p)}$-local
stable homotopy category has been given by Franke~\cite{franke}.  The
interested reader may wish to consult~\cite{roitzheim, ganter}.  In a
different direction, we note the further work of Bousfield, building a
unified version of \Ktheory\ in order to combine information from
different primes~\cite{Bousfield90}.  

\medskip

Throughout this paper $p$ will denote an odd prime.

\medskip

We would like to thank Peter Kropholler for pointing out a result
which we need in the proof of Proposition~\ref{globalfree}.  The third
author acknowledges the support of a Scheme~4 grant from the London
Mathematical Society.

\section{Discrete $A$-modules}
\label{DiscreteModules}

Let $p$ be an odd prime and let $A=K_{(p)}^{0}(K_{(p)})$ be the ring
of degree zero stable operations in \plocal\ complex \Ktheory.  This
ring can be described as follows; see~\cite{CCW05}.  Choose an
integer~$q$ that is primitive modulo~$p^{2}$, and let $\Psiq\in A$ be
the corresponding Adams operation.  Let
\begin{equation}\label{qidefn}
    q_{i}
    =
    q^{(-1)^{i}\floor{i/2}},
\end{equation}
and define polynomials $\Theta_{n}(X)$, for each integer $n\ge 0$, by
$\Theta_{n}(X)=\prod_{i=1}^{n}(X-q_{i})$.  Then let the operation
$\Phi_{n}\in A$ be given by $\Phi_{n}=\Theta_{n}(\Psiq)$.  For
example, $\Phi_{4} = (\Psiq-1)(\Psiq-q)(\Psiq-q^{-1})(\Psiq-q^{2})$.
These operations have been chosen so that any infinite sum 
$\sum_{n\ge 0}a_{n}\Phi_{n}$, with coefficients $a_{n}$ in the
\plocal\ integers~$\Zpl$, converges.  The following theorem says that
any operation can be written uniquely in this form.

\begin{thm}\cite[Theorem~6.2]{CCW05}
The elements of $A$ can be expressed uniquely as infinite sums
$\sum_{n\ge 0}a_{n}\Phi_{n}$, where $a_{n}\in\Zpl$.
\qed
\end{thm}

For each $m\geq 0$, we define
$$
A_{m}
=
\biggl\lbrace\,
\sum_{n\ge m}a_{n}\Phi_{n}:a_{n}\in\Zpl
\biggr\rbrace
\subseteq
A,
$$
so that $A_{m}$ is the ideal of operations which annihilate the
coefficient groups $\pi_{2i}(K_{(p)})$ for $-m/2<i<(m+1)/2$, and thus
it does not depend on the choice of primitive element~$q$.  We obtain
a decreasing filtration
$$
A
=
A_{0}
\supset
A_{1}
\supset\dots\supset
A_{m} 
\supset
A_{m+1}
\supset
\dotsb.
$$ 

We use this filtration to give the ring~$A$ a topology in the standard
way (see, for example, Chapter~9 of~\cite{Northcott68}): the open sets
are unions of sets of the form $y+A_{m}$, where $y\in A$ and $m\ge 0$.
We note that $A$ is complete with respect to this topology.  Indeed,
another way to view the topological ring $A$ is as the completion of
the polynomial ring $\Zpl[\Psiq]$ with respect to the filtration by
the principal ideals $\Phi_{m}\Zpl[\Psiq]=\Zpl[\Psiq]\cap A_{m}$.
(By way of warning, this filtration is not multiplicative in the sense
of~\cite{Northcott68}, and, since $A$ is not
Noetherian~\cite[Theorem~6.10]{CCW05}, it is not the completion of a
polynomial ring with respect to \emph{any} multiplicative filtration.)

\begin{defn}\label{discrete}
    A \emph{discrete $A$-module} $M$ is an $A$-module such that the
    action map
    $$
    A\times M\to M
    $$
    is continuous with respect to the discrete topology on $M$ and the
    resulting product topology on $A\times M$.
\end{defn}

In practice we use the following criterion to recognise discrete 
modules.

\begin{lemma}\label{discretecriterion}
    An $A$-module $M$ is discrete if and only if for each $x\in M$,
    there is some~$n$ such that $A_{n}x=0$.
\end{lemma}

\begin{proof}
    Fixing $x\in M$, the map sending $\alpha\in A$ to 
    $(\alpha,x)\in A\times M$ is continuous.  Thus if $M$ is 
    discrete, the map $\alpha\mapsto\alpha x\in M$ is continuous, 
    which implies that its kernel contains $A_{n}$ for some~$n$.
    
    Suppose now that for each $x\in M$, there exists~$n$ such that
    $A_{n}x=0$.  If $x,y\in M$ are such that $\alpha x=y$, and 
    $A_{n}x=0$, then $(\alpha+A_{n})x=y$ so that 
    $(\alpha+A_{n})\times\lbrace x\rbrace$ is an open neighbourhood 
    of $(\alpha,x)$ in the preimage under the action map 
    of~$\lbrace y\rbrace$.  This shows that $M$ is a discrete module.
\end{proof}

\begin{remark}
    For $n\ge 1$, the principal ideal $\Phi_{n}A$ is \emph{strictly}
    contained in $A_{n}$; see Proposition~\ref{An/PhinA}.  Thus, if an
    $A$-module $M$ has the property that for each $x\in M$ there is
    some $n\ge 0$ with $\Phi_{n}x=0$, then it does not follow that $M$
    is a discrete $A$-module.

    For example, the $A$-module $A/\Phi_{1}A$ is not discrete,
    although every element is annihilated by $\Phi_{1}$.  For
    suppose the element $1+\Phi_{1}A\in A/\Phi_{1}A$ were annihilated
    by the ideal $A_{n}$, i.e., $A_{n}\subseteq\Phi_{1}A$.  Then
    $A/\Phi_{1}A$ would be finitely generated over~\Zpl, being a
    quotient of $A/A_{n}$ which has rank~$n$ over~\Zpl.  However,
    $A/\Phi_{1}A$ contains the submodule $A_{1}/\Phi_{1}A$ which is
    isomorphic to~$\Zp/\Zpl$, where $\Zp$ denotes the \padic\
    integers.  The proof of this is analogous to that of the
    corresponding result for connective \Ktheory\ which is given in
    \cite[Corollary~3.7]{CCW05}.

    We will see in Section~\ref{localfg} that if $Ax$ is finitely
    generated over~\Zpl, then $\Phi_{n}x=0$ does imply $A_{n}x=0$.
\end{remark}

The motivating example for the concept of a discrete $A$-module is the
\plocal\ $K$-homology of a spectrum.

\begin{prop}\label{spectraarecontinuous}
    The degree zero \plocal\ $K$-homology $K_{0}(X;\Zpl)$ of a
    spectrum $X$ is a discrete $A$-module with the action map given by
    $$
    A\tensor K_{0}(X;\Zpl)\to 
    A\tensor\KzKp\tensor K_{0}(X;\Zpl)\to
    K_{0}(X;\Zpl),
    $$
    in which \KzKp\ denotes the ring of degree zero cooperations in
    \plocal\ \Ktheory, the first map arises from the coaction map, and
    the second comes from the Kronecker pairing.
\end{prop}

\begin{proof}
    Consider the Kronecker pairing 
    $\langle\_,\_\rangle:A\tensor\KzKp\to\Zpl$.  We claim that for
    each element $f\in\KzKp$, there is some~$m$ such that 
    $\langle A_{m},f\rangle=0$.  Indeed this is clear, since
    $A=\Homz(\KzKp,\Zpl)$, with the $\Phi_{n}$ being the dual basis to
    a particular \Zpl-basis of $\KzKp$; see~\cite[Theorem~6.2]{CCW05}.

    Now let $X$ be a spectrum, and let $x\in K_{0}(X;\Zpl)$.  Under
    the coaction $K_{0}(X;\Zpl)\to\KzKp\tensor K_{0}(X;\Zpl)$, the
    image of~$x$ is a finite sum $\sum_{i=1}^{k}f_{i}\tensor x_{i}$
    for some elements $f_{i}\in\KzKp$ and $x_{i}\in K_{0}(X;\Zpl)$.
    By the preceding paragraph, for each $f_{i}$, there is
    some~$m_{i}$ such that $\langle A_{m_{i}},f_{i}\rangle=0$.  Hence,
    if $m=\max\{m_{1},\dots,m_{k}\}$, we have 
    $\langle A_{m},f_{i}\rangle=0$ for all $i$.  Thus
    $A_{m}x=\sum_{i=1}^{k}\langle A_{m},f_{i}\rangle x_{i}=0$.
\end{proof}

If an $A$-module is a discrete $A$-module, then the $A$-action is
determined by the action of $\Phi_{1}$, or, equivalently, by the action
of $\Psiq=1+\Phi_{1}$.

\begin{lemma}\label{determinedbyphi1}
    If $M$ and $N$ are discrete $A$-modules, and $f:M\to N$ is a
    \Zpl-homomorphism which commutes with the action of~\Psiq,
    then $f$ is a homomorphism of $A$-modules.
\end{lemma}

\begin{proof}
    Since $\Phi_{i}=\Theta_{i}(\Psiq)$, the map $f$ must also commute
    with~$\Phi_{i}$ and hence with all finite linear combinations
    $\sum_{i=0}^{k}a_{i}\Phi_{i}$.  Now let 
    $\alpha=\sum_{i\ge 0}a_{i}\Phi_{i}\in A$ and $x\in M$.  Since $M$
    and~$N$ are discrete $A$-modules, for some $n$, we have $A_{n}x=0$ and
    $A_{n}f(x)=0$.  Then 
    $$
    f\biggl(\sum_{i\ge 0}a_{i}\Phi_{i}x\biggr) 
    = 
    f\biggl(\sum_{i=0}^{n-1}a_{i}\Phi_{i}x\biggr) 
    = 
    \sum_{i=0}^{n-1}a_{i}\Phi_{i}f(x) 
    = 
    \sum_{i\ge 0}a_{i}\Phi_{i}f(x).
    $$ 
    Hence $f$ commutes with all elements of $A$, i.e., it is an
    $A$-ho\-mo\-mor\-phism.
\end{proof}

We end this section by making some remarks about the category of
discrete $A$-modules.

\begin{defn}
    We define \emph{the category of discrete $A$-modules} \DA\ to
    be the full subcategory of the category of $A$-modules whose
    objects are discrete $A$-modules.
\end{defn}

\begin{lemma}\label{submodulesandquotients}
    The category \DA\ is closed under submodules and quotients.
\end{lemma}

\begin{proof}
    It is clear that a submodule of a discrete $A$-module is discrete.

    Let $M$ be a discrete $A$-module and suppose that $N$ is an
    $A$-submodule of~$M$.  Consider the quotient $A$-module $M/N$.
    For $x\in M$, there is some~$n$ such that $A_{n}x=0$.  Then
    $A_{n}(x+N)=0$, and so $M/N$ is a discrete $A$-module.
\end{proof}

\begin{lemma}\label{directsums}
    The category \DA\ has arbitrary direct sums. 
\end{lemma}

\begin{proof}
    Let $M_{i}$ be a discrete $A$-module for each $i$ in some indexing set
    $\mathcal{I}$.  Each element of $\oplus_{i\in\mathcal{I}}M_{i}$ is
    a finite sum $x_{1}+\dots+x_{k}$, where each $x_{j}$ belongs to
    some $M_{i}$.  For each $j$, there is some $n_{j}$ such that
    $A_{n_{j}}x_{j}=0$.  Hence $A_{n}(x_{1}+\dots+x_{k})=0$ if
    $n\ge\max\{n_{1},\dots,n_{k}\}$.  So the direct sum
    $\oplus_{i\in\mathcal{I}}M_{i}$ is discrete.
\end{proof}

\begin{cor}\label{abelian}
    \DA\ is a cocomplete abelian category. 
\end{cor}

\begin{proof}
    That \DA\ is abelian follows directly from
    Lemmas~\ref{submodulesandquotients} and~\ref{directsums} and the
    fact that $A$-modules form an abelian category.  An abelian
    category with arbitrary direct sums is cocomplete; see, for
    example,~\cite[Proposition~2.6.8]{weibel}.
\end{proof}

When combined with Theorem~\ref{isomorphism} below, the following
result corresponds to~10.5 of~\cite{Bousfield85}.

\begin{prop}
    \DA\ is isomorphic to the category of \KzKp-comodules.
\end{prop}

\begin{proof}
    Let $G_{n}(w)=q^{n\floor{n/2}}F_{n}(w)\in\KzKp$, where $F_{n}(w)$ 
    is as defined in the proof of Theorem~6.2 in~\cite{CCW05}.  That 
    proof shows that the $\Phi_{n}$ are dual to the $G_{n}(w)$, i.e., 
    $\langle\Phi_{n},G_{i}(w)\rangle=\delta_{ni}$.
    
    If $M$ is a discrete $A$-module, define 
    $\phi_{M}:M\to\KzKp\tensor M$ by 
    $\phi_{M}(x)=\sum_{n\ge 0}G_{n}(w)\tensor\Phi_{n}x$, in which only
    finitely many terms are non-zero because $M$ is discrete.  It is
    easily checked that this makes $M$ into a \KzKp-comodule and that
    an $A$-module homomorphism between discrete $A$-modules is also a
    \KzKp-comodule homomorphism.
     
    If $M$ is a \KzKp-comodule with coaction 
    $\phi_{M}:M\to\KzKp\tensor M$, then, just as in the proof of 
    Proposition~\ref{spectraarecontinuous}, $M$ is a discrete 
    $A$-module with action given by
    $$
    A\tensor M\xrightarrow{\;1\tensor\phi_{M}\;}
    A\tensor\KzKp\tensor M\xrightarrow{\;\langle\_,\_\rangle\tensor 1\;}
    M.
    $$
    It is clear that this construction is functorial.
    
    It is routine to check that the two constructions are mutually 
    inverse.    
\end{proof}

\section{Locally Finitely Generated $A$-modules}
\label{localfg}

\begin{defn}
    An $A$-module $M$ is \emph{locally finitely generated} if, for
    every $x\in M$, the submodule $Ax$ is finitely generated
    over~$\Zpl$.
\end{defn}

This section is devoted to proving the following theorem.

\begin{thm}\label{lfg}
    An $A$-module $M$ is discrete if and only if it is locally
    finitely generated.
\end{thm}

We need to prove first a number of preliminary results.

\begin{prop}\label{An/PhinA}
    For all $n\ge 1$, the quotient $A_{n}/\Phi_{n}A$ is a 
    rational vector space.
\end{prop}

\begin{proof}
    The result will follow if we can show that $A_{n}/\Phi_{n}A$ is
    divisible and torsion-free.  We thus need to show:
    \begin{enumerate}
        \item\label{pdiv} 
        For any $\alpha\in A_{n}$, there exists $\beta\in A$ such that
        $\Phi_{n}\beta-\alpha\in pA_{n}$;
    
        \item\label{ptf}
        if $\alpha\in A_{n}$ and $p\alpha\in\Phi_{n}A$, then
        $\alpha\in\Phi_{n}A$.
    \end{enumerate}
    Recall from Proposition~6.8 of~\cite{CCW05} that
    $$
    \Phi_{n}\Phi_{j}
    =
    \!\sum_{k=\max(j,n)}^{j+n}\!c_{j,n}^{k}\Phi_{k}
    $$
    for certain coefficients $c_{j,n}^{k}\in\Zpl$ (denoted 
    $A_{j,n}^{k}$ in~\cite{CCW05}).
    
    Suppose $\alpha=\sum_{k\ge n}a_{k}\Phi_{k}\in A_{n}$, then
    $\beta=\sum_{j\ge 0}b_{j}\Phi_{j}\in A$ will satisfy  
    $\Phi_{n}\beta-\alpha\in pA_{n}$ if and only if the  
    congruences
    \begin{equation}\label{abcongs}
        \sum_{j=k-n}^{k}c_{j,n}^{k}b_{j}
        \equiv
        a_{k}
        \mod{p}
        \quad
        (k\ge n)
    \end{equation}
    can be solved for $(b_{j})_{j\ge 0}$.  We will verify~\eqref{pdiv}
    by showing that these congruences \emph{can} always be solved,
    and~\eqref{ptf} by showing that the solution is unique modulo~$p$.
    
    It follows from Propositions~A.2 and~A.4 of~\cite{CCW05} that, 
    with $q_{k}$ as defined in~\eqref{qidefn},
    \begin{equation}\label{crecurs}
        c_{j,n}^{k}
        =
        (q_{k+1}-q_{n})c_{j,n-1}^{k}+c_{j,n-1}^{k-1},
    \end{equation}
    where $c_{j,n}^{k}=0$ unless $j,n\le k\le j+n$, and
    $c_{j,n}^{k}=1$ if $k=j+n$.  We claim that for any integer 
    $s\ge 1$,
    \begin{equation}\label{ccong}
        c_{j,n}^{k}
        \equiv
        0\mod{p}
        \quad
        \text{for $(2p-2)s\le j\le k<(2p-2)s+n$.}
    \end{equation}
    This follows by induction on~$n$ from~\eqref{crecurs} and the
    periodicity of the sequence $(q_{i})_{i\ge 1}$ modulo~$p$:
    $q_{i}\equiv q_{i+2p-2}\mod{p}$.
    
    It follows from~\eqref{ccong} that for $k=(2p-2)s+n-1$ the
    congruence~\eqref{abcongs} is 
    $b_{(2p-2)s-1}\equiv a_{(2p-2)s+n-1}\mod{p}$, and that for
    $k<(2p-2)s+n-1$ the congruences have the form
    $$
    b_{k-n}
    +
    \bigl(\text{terms involving $b_{j}$ for $k-n<j<(2p-2)s$}\bigr)
    \equiv
    a_{k}
    \mod{p}.
    $$
    It is now clear that the congruences~\eqref{abcongs} have a unique 
    solution modulo $p$ for the~$b_{j}$.
\end{proof}

\begin{cor}\label{phitoA}
    If $M$ is a locally finitely generated $A$-module and $x\in M$
    satisfies $\Phi_{n}x=0$ for some $n\ge 0$, then $A_{n}x=0$.
\end{cor}

\begin{proof}
    By hypothesis the map $A_{n}\to Ax$ sending $\alpha$ to~$\alpha x$ 
    factors through the \Q-module $A_{n}/\Phi_{n}A$. Hence, since $Ax$ 
    is a finitely generated \Zpl-module, the map must be zero.
\end{proof}

We now need to show that, for an element $x$ in a locally finitely
generated $A$-module $M$, we have $\Phi_{n}x=0$ for some $n\ge 0$.
The next results establish this, first for \Zpl-torsion modules, then
for $\Zpl$-free modules and finally in the general case.

\begin{prop}\label{finiteA-mods}
    If $M$ is a finite $A$-module, then $\Phi_{n}M=0$ for some~$n$.
\end{prop}

\begin{proof}
    Let $x\in M$ be non-zero, and define
    $$
    I
    =
    \set{f(X)\in\Zpl[X]:f(\Psiq)x=0},
    $$
    which is clearly an ideal of~$\Zpl[X]$.
    
    Let $f(X)$ be any element of~$I$ such that its reduction
    $\fbar(X)\in\Fp[X]$ is not zero.  Since the elements
    $\Psi^{q^{r}}\!x$ for $r\ge 0$ cannot be distinct, we may, for
    example, take an $f(X)$ of the form $X^{r_{1}}-X^{r_{2}}$ with
    $r_{1}\ne r_{2}$.
    
    Suppose that
    $$
    \fbar(X)
    =
    \gbar(X)\prod_{k=1}^{p-1}(X-k)^{e_{k}}
    \qquad
    (e_{k}\ge 0),
    $$
    where $\gbar(X)$ has no roots in~\Fpu, and thus we can write
    $$
    f(X)
    =
    g(X)\prod_{k=1}^{p-1}(X-k)^{e_{k}}+ph(X),
    $$
    for some $g(X),h(X)\in\Zpl[X]$, where $g(k)\in\Zplu$ for 
    $k=1,2,\dots,p-1$.
    
    We recall now from \cite{Johnson87} and \cite[\S6]{CCW05} that if
    $\phi(X)\in\Zpl[X]$, the element $\phi(\Psiq)$ is a unit in~$A$ if
    and only if $\phi(q_{i})$ is a unit in~\Zpl\ for all $i\ge 1$.
    Since the $q_{i}$ take the values $1,2,\dots,p-1$ modulo~$p$, it
    follows that $g(\Psiq)$ is a unit in~$A$.  Since $f(\Psiq)$ cannot
    be a unit, it follows that $e_{k}>0$ for at least one value
    of~$k$.
    
    We have
    $$
    g(\Psiq)\prod_{k=1}^{p-1}(\Psiq-k)^{e_{k}}x
    =
    -ph(\Psiq)x,
    $$
    and thus
    $$
    \prod_{k=1}^{p-1}(\Psiq-k)^{e_{k}}x
    =
    p\alpha x,
    $$
    where $\alpha=-g(\Psiq)^{-1}h(\Psiq)\in A$.
    
    As $M$ is finite, there is some~$s$ such that $p^{s}x=0$, so that
    $$
    \prod_{k=1}^{p-1}(\Psiq-k)^{e_{k}s}x
    =
    p^{s}\alpha^{s}x
    =
    0.
    $$
    
    But it is clear that $\prod_{k=1}^{p-1}(X-k)^{e_{k}s}$ is a 
    factor modulo~$p^{s}$ of $\Theta_{n}(X)$ for sufficiently 
    large~$n$, and thus $\Phi_{n}x=0$ for such~$n$.
    
    By choosing the maximum such~$n$ over all non-zero $x\in M$, we 
    have $\Phi_{n}M=0$.
\end{proof}

\begin{prop}\label{globalfree}
    If $M$ is an $A$-module which is free of finite rank over~\Zpl,
    then $\Phi_{n}M=0$ for some~$n$.
\end{prop}

\begin{proof}
    Let $x\in M$.  The action map $\eta_{x}:A\to M$ given by
    $\eta_{x}(\alpha)=\alpha x$ is a homomorphism of $A$-modules.  In
    particular, it is a homomorphism of abelian groups.
    By~\cite[Theorem~95.3]{fuchsvol2}, the target is a slender group,
    so there is some~$m$ such that $\eta_{x}(\Phi_{m})=0$, i.e., 
    $\Phi_{m}x=0$.

    If $x_{1},\dots,x_{r}$ is a \Zpl-basis of $M$ and
    $\Phi_{m_{i}}x_{i}=0$ for $i=1,\dots,r$, then $\Phi_{n}M=0$ where
    $n=\max\{m_{i}:1\le i\le r\}$.
\end{proof}

Since $A$ is not an integral domain, we must exercise some care with 
quotients.  However, if $n>m$ the polynomial $\Theta_{m}(X)$ is a 
factor of $\Theta_{n}(X)$, so we may let $\Phi_{n}/\Phi_{m}$ denote 
the value of the polynomial $\Theta_{n}(X)/\Theta_{m}(X)$ at~\Psiq.

\begin{lemma}\label{PhiQuotientsModp}
    If $n>m$ and $n-m$ is divisible by $2p-2$, then 
    $$
    \frac{\Phi_{n}}{\Phi_{m}} 
    \equiv\Phi_{n-m}
    \mod{p^{1+\nu_{p}(n-m)}}.
    $$
\end{lemma}

\begin{proof}
    It is easy to verify that whenever $n>m>0$,
    $$
    \frac{\Phi_{n}}{\Phi_{m}} 
    =
    \frac{\Phi_{n-1}}{\Phi_{m-1}}+(q_{m}-q_{n})\frac{\Phi_{n-1}}{\Phi_{m}}.
    $$
    Expanding the first term on the right in the same way, and
    repeating this process, leads to the equation
    $$
    \frac{\Phi_{n}}{\Phi_{m}} 
    = 
    \Phi_{n-m}+\sum_{i=0}^{m-1}(q_{m-i}-q_{n-i})\frac{\Phi_{n-i-1}}{\Phi_{m-i}}.
    $$
    If $n-m = 2k$, then $q_{m-i}-q_{n-i}$ is divisible by $q^{k}-1$
    for all~$i$.  Moreover, if $k$ is divisible by $(p-1)p^{r-1}$,
    then $q^{k}-1\equiv 0\mod{p^{r}}$.
\end{proof}

\begin{proof}[Proof of Theorem~\ref{lfg}]
    First suppose that $M$ is a discrete $A$-module.  Then, for 
    $x\in M$, there is some $n$ such that $Ax=(A/A_{n})x$.  But
    $A/A_{n}$ is free of finite rank over~\Zpl, so $M$ is locally
    finitely generated.

    Now suppose that $M$ is locally finitely generated.  By
    Corollary~\ref{phitoA}, it is enough to show that for each 
    $x\in M$ there is some $n$ such that $\Phi_{n}x=0$.

    Let $x\in M$, and write $N=Ax$.  By hypothesis, $N$ is finitely
    generated over~\Zpl.  Let $T\subseteq N$ be the $A$-submodule of
    $N$ consisting of the \Zpl-torsion elements.  So there is some $s$
    such that $p^{s}T=0$.  The quotient $N/T$ is an $A$-module which is
    free of finite rank over~\Zpl.  By
    Proposition~\ref{globalfree}, there is some~$k$ such that
    $\Phi_{k}(N/T)=0$, and so $\Phi_{k}x\in T$.  Then, by
    Proposition~\ref{finiteA-mods}, there is some~$r$ such that
    $\Phi_{r}\Phi_{k}x=0$.

    By increasing~$r$ if necessary, we may arrange that $r$ is
    divisible by $(2p-2)p^{s}$, so that, by
    Lemma~\ref{PhiQuotientsModp}, we have
    $$
    \frac{\Phi_{r+k}}{\Phi_{k}}
    =
    \Phi_{r}+p^{s}\theta, 
    $$
    for some $\theta\in A$.  Hence 
    $$
    \Phi_{r+k}x 
    =
    \frac{\Phi_{r+k}}{\Phi_{k}}\Phi_{k}x 
    = 
    \Phi_{r}\Phi_{k}x+p^{s}\theta\Phi_{k}x
    =
    0+\theta(p^{s}\Phi_{k}x) 
    =
    0.
    $$
\end{proof}

\section{Bousfield's Category of $K$-theory modules}
\label{bousfieldscategory}

In this section we relate the category \DA\ of discrete $A$-modules
to a category considered by Bousfield in his work on the
$K_{(p)}$-local stable homotopy category.

Let $R=\Zpl[\Zplu]$ be the group-ring of the multiplicative group of
units in $\Zpl$, with coefficients in $\Zpl$ itself.  For clarity, as
well as to reflect the topological applications, we write 
$\Psi^{j}\in R$ for the element $j\in\Zplu$.  Hence elements of $R$
are finite $\Zpl$-linear combinations of the~$\Psi^{j}$.

\begin{defn}\cite{Bousfield85}\label{bousfieldDefn}
    \emph{Bousfield's category} $\A(p)$ is the full subcategory of the
    category of $R$-modules whose objects $M$ satisfy the following
    conditions.  For each $x\in M$,
    \begin{itemize}
        \item[(a)] the submodule $Rx\subseteq M$ is finitely
        generated over~\Zpl,

        \item[(b)] for each $j\in\Zplu$, $\Psi^{j}$ acts on
        $Rx\tensor\Q$ by a diagonalisable matrix whose eigenvalues are
        integer powers of~$j$,

        \item[(c)] for each $m\ge 1$, the action of \Zplu\ on
        $Rx/p^{m}Rx$ factors through the quotient homomorphism
        $\Zplu\to(\Z/p^{k}\Z)^{\times}$ for sufficiently large~$k$.
    \end{itemize}
    We call the objects of this category \emph{Bousfield modules}.
\end{defn}

If $M$ is a Bousfield module which is finitely generated over $\Zpl$,
then condition (a) holds automatically, and the rational
diagonalisability and \padic\ continuity conditions of (b) and (c)
hold globally for $M$, as well as for each submodule $Rx\subseteq M$.

\medskip

Note that $R\subset A$, so an $A$-module can be considered as an
$R$-module by restricting the action.  (An explicit formula expressing
each $\Psi^{j}$, for $j\in\Zplu$, in terms of the $\Phi_{n}$ is given
in~\cite[Proposition~6.6]{CCW05}.)

\begin{thm}\label{continuousimpliesbousfield}
    If $M$ is a discrete $A$-module, then $M$ is a Bousfield module
    with the $R$-action given by the inclusion $R\subset A$.
\end{thm}

\begin{proof}
    Suppose $A_{n}x=0$, then it is clear that $Rx=Ax=(A/A_{n})x$,
    which is finitely generated over~\Zpl.  The matrix representing
    \Psiq\ on $Rx\tensor\Q$ is annihilated by the polynomial
    $\Theta_{n}(X)$, and thus its minimal polynomial is a factor of
    $\Theta_{n}(X)$.  Since $\Theta_{n}(X)$ has distinct rational
    roots, the matrix can be diagonalised over~\Q. The eigenvalues are
    roots of $\Theta_{n}(X)$, which are integer powers of~$q$.
    
    Proposition~6.6 of~\cite{CCW05} shows that
    $$
    \Psi^{j}
    =
    \sum_{i\ge 0}g_{i}(j)\Phi_{i},
    $$
    if $j\in\Zplu$, where $g_{i}(w)$ is a certain Laurent polynomial
    (given explicitly in~\cite{CCW05}) satisfying $g_{i}(j)\in\Zpl$
    for all $j\in\Zplu$.  Since $A_{n}x=0$, this shows that $\Psi^{j}$
    acts on~$Rx$ as the polynomial in~\Psiq
    $$
    P_{j}(\Psiq)
    =
    \sum_{i=0}^{n-1}g_{i}(j)\Theta_{i}(\Psiq).
    $$
    It follows that the eigenvalues of the matrix of the action
    of~$\Psi^{j}$ on $Rx\tensor\Q$ are $P_{j}(q^{r})$, where $q^{r}$
    is an eigenvalue of the action of~\Psiq.  By considering the
    action on the coefficient group $\pi_{2r}(K_{(p)})$, we see that
    $P_{j}(q^{r})=j^{r}$.  Hence condition~(b) of
    Definition~\ref{bousfieldDefn} is satisfied.

    The Laurent polynomials $g_{i}(j)$ are uniformly \padically\
    continuous functions of $j\in\Zplu$, i.e., for each $m\ge 1$ there
    is an integer $K_{i}$ such that $g_{i}(j)\equiv
    g_{i}(j+p^{k}a)\mod{p^{m}}$ whenever $k\ge K_{i}$.  Hence
    $\Psi^{j}x\equiv\Psi^{j+p^{k}a}x\mod{p^{m}}$ for
    $k\ge\max\{K_{0},K_{1},\dots,K_{n-1}\}$.  This shows that
    condition~(c) of Definition~\ref{bousfieldDefn} holds.
\end{proof}

Since each $\Phi_{n}$ is a polynomial in~\Psiq, any finite linear
combination of the $\Phi_{n}$ can be considered as an element of~$R$.
In order to show that a Bousfield module can be given the structure of
an $A$-module, we need to specify how an infinite sum 
$\sum_{n\ge 0}a_{n}\Phi_{n}$ acts.  We need a preliminary lemma.

\begin{lemma}\label{algcondmodp}
    Let $M$ be a Bousfield module and $x\in M$.  There is some $k\ge
    1$ such that $\Phi_{n}x\equiv 0\mod{p}$ for all $n\ge p^{k}(p-1)$.
\end{lemma}

\begin{proof}
    Proposition~6.5 of~\cite{CCW05} gives an explicit formula for the
    expansion of $\Phi_{p^{k}(p-1)}$ as a finite $\Zpl$-linear
    combination of the $\Psi^{q^{j}}$.  The coefficient of $\Psi^{q^{j}}$
    has as a factor the $q$-binomial coefficient $\gp{p^{k}(p-1)}{j}$,
    which is divisible by~$p$ for $0<j<p^{k}(p-1)$.  Thus
    \begin{align*}
        \Phi_{p^{k}(p-1)} 
        &\equiv 
        \Psi^{q^{p^{k}(p-1)}}+q^{p^{k}(p-1)/2}\mod{p}\\
        &\equiv
        \Psi^{q^{p^{k}(p-1)}}-1\mod{p}.
    \end{align*}
    Now condition~(c) of Definition~\ref{bousfieldDefn} ensures that
    $(\Psi^{q^{p^{k}(p-1)}}-1)x\equiv 0\mod{p}$ for sufficiently
    large~$k$.  Thus $\Phi_{p^{k}(p-1)}x\equiv 0\mod{p}$, and
    consequently $\Phi_{n}x\equiv 0\mod{p}$ for all $n\ge p^{k}(p-1)$.
\end{proof}

\begin{thm}\label{bousfieldimpliescontinuous}
    If the $R$-module $M$ is a Bousfield module, then the $R$-action
    extends uniquely to an $A$-action in such a way as to make $M$ a
    discrete $A$-module.
\end{thm}

\begin{proof}
    Let $x\in M$.  By condition~(b) of Definition~\ref{bousfieldDefn},
    the minimal polynomial of the action of \Psiq\ on
    $Rx\tensor\Q$ has the form $\prod_{i=1}^{t}(X-q^{k_{i}})$, where
    $k_{i}\in\Z$.  For sufficiently large~$m$, this polynomial is a
    factor of $\Theta_{m}(X)$, so that $\Phi_{m}x=0$ in $Rx\tensor\Q$.
    This means that $\Phi_{m}x\in T$, the \Zpl-torsion submodule
    of~$Rx$.  As $T$ is finitely generated, there is some exponent~$e$
    such that $p^{e}T=0$.
    
    Let $k$ be as in Lemma~\ref{algcondmodp}, and let 
    $n=(2p-2)p^{r}$, where $r\ge k$.  Then $\Phi_{n}x\equiv 0\mod{p}$, 
    and so $\Phi_{n}^{\ell}x\equiv 0\mod{p^{\ell}}$ for any $\ell\ge 1$.
    
    On the other hand, iterating Lemma~\ref{PhiQuotientsModp} shows
    that $\Phi_{2^{s}n}\equiv\Phi_{n}^{2^{s}}\mod{p^{1+r}}$.  It is
    thus clear that by choosing $r$ and~$s$ sufficiently large we may
    ensure firstly that $\Phi_{2^{s}n}x\in T$ and then that
    $\Phi_{2^{s}n}x\equiv 0\mod{p^{e}}$, which means that
    $\Phi_{2^{s}n}x=0$.  
    
    There is thus some integer $N$ such that $\Phi_{N}x=0$.  Given 
    $\alpha=\sum_{k\ge 0}a_{k}\Phi_{k}\in A$, let 
    $\alpha x=\sum_{k=0}^{N-1}a_{k}\Phi_{k}x$.  It is clear that this 
    gives~$M$ the structure of a discrete $A$-module.  The uniqueness 
    of this structure follows from Lemma~\ref{determinedbyphi1}.
\end{proof}

\begin{remark}
    We have chosen to give a direct algebraic proof.  This can, of
    course, be replaced by appealing to Bousfield's topological
    results: the non-split analogue of
    ~\cite[Proposition~8.7]{Bousfield85} shows that any object in
    $\A(p)$ can be expressed as $K_{0}(X;\Zpl)$ for some spectrum~$X$.
    Proposition~\ref{spectraarecontinuous} shows that $K_{0}(X;\Zpl)$
    is a discrete $A$-module.
\end{remark}

\medskip

Assembling the main results of this section, we have proved the
following.

\begin{thm}\label{isomorphism}
    Bousfield's category $\A(p)$ is isomorphic to the category~\DA\
    of discrete $A$-modules.  
    \qed
\end{thm}

\section{Cofree Objects and Injective Resolutions}
\label{cofreeobjs}

To illustrate the utility of our point of view, we consider cofree
objects.  Bousfield showed the existence of a right adjoint functor
$U$ to the forgetful functor from his category $\A{(p)}$ to the
category of $\Zpl$-modules.  To do this, he had to give different
descriptions in two special cases and then deduce the existence of
such a functor in the general case without constructing it explicitly.

In our context, the functor $U$ is just a continuous $\Hom$ functor,
right adjoint to the forgetful functor from \DA\ to $\Zpl$-modules.
So, not only do we not need to treat separate cases, but our uniform
description is conceptually simple and fits into a standard framework
for module categories.

If $S$ is a unital, commutative algebra over a commutative ring $R$
with~$1$, then there is a right adjoint to the forgetful functor from
the category of $S$-modules to the category of $R$-modules, given by
$\Hom_{R}(S,\_)$; see~\cite[Lemma~2.3.8]{weibel}.  For an $R$-module $L$,
$\Hom_R(S,L)$ is an $S$-module via $(sf)(t)=f(ts)$ for $s,t\in S$.

We will need to modify this construction since $\Homz(A,L)$
need not be a \emph{discrete} $A$-module.

\begin{ex}\label{Qex}
    The $A$-module $\Homz(A,\Q)$ is not discrete.  To see this, note
    that, since $\Q$ is not slender~\cite[Section 94]{fuchsvol2},
    there exists $f\in\Homz(A,\Q)$ such that $f(\Phi_{n})\ne 0$ for
    all~$n$.  Then $(\Phi_{n}f)(1)=f(\Phi_{n})\ne 0$ for all $n$, and
    so there is no $n$ such that $A_{n}f=0$.

    Notice that we can make \Q\ into a discrete $A$-module via
    $A\xrightarrow{\epsilon}\Zpl\hookrightarrow\Q$, where
    $A\xrightarrow{\epsilon}\Zpl$ is the augmentation given by
    $\epsilon\bigl(\sum_{k\ge 0}a_{k}\Phi_{k}\bigr)=a_{0}$.  Thus
    $\Homz(A,L)$ need not be a discrete $A$-module, even when $L$ is
    such a module.
\end{ex}

\begin{defn}
    If $L$ is a \Zpl-module, let $\Homc(A,L)$ denote the 
    $A$-submodule of $\Homz(A,L)$ consisting of the homomorphisms which 
    are continuous with respect to the filtration topology on~$A$ and 
    the discrete topology on~$L$.  
\end{defn}

Note that a homomorphism $A\to L$ is continuous if and only if its
kernel contains $A_{n}$ for some $n$.  Example~\ref{Qex} shows that
$\Homc(A,L)$ may be a proper submodule of~$\Homz(A,L)$ even if $L$ is
a discrete $A$-module.

\begin{prop}\label{cofreeconstruction}
    The functor $U$ from $\Zpl$-modules to the category \DA\ of
    discrete $A$-modules given by $UL=\Homc(A,L)$ is
    right adjoint to the forgetful functor.
\end{prop}

\begin{proof}
    For any $A$-module $N$, we define the discrete heart of $N$ to be
    the largest $A$-submodule which is a discrete $A$-module:
    $$
    N^{\text{disc}}
    =
    \set{x\in N:\text{$A_{n}x=0$ for some $n$}}.
    $$
    It is easy to show that the `discrete heart functor' is right
    adjoint to the forgetful functor from discrete $A$-modules to
    $A$-modules, hence $\Homz(A,\_)^{\operatorname{disc}}$ is right
    adjoint to the forgetful functor from discrete $A$-modules to
    \Zpl-modules.  The proof is completed by observing that
    $\Homz(A,L)^{\operatorname{disc}}=\Homc(A,L)$ for any \Zpl-module
    $L$.
\end{proof}

\begin{prop}\label{coopsuniversal}
    For any \Zpl-module $L$, there is a natural isomorphism of
    $A$-modules $UL\cong K_{0}(K)_{(p)}\tensor L$, where
    $K_{0}(K)_{(p)}\tensor L$ is an $A$-module via the $A$-module
    structure on~$K_{0}(K)_{(p)}$.
\end{prop}

\begin{proof}
    The map $\KzKp\to\Homc(A,\Zpl)$ which sends $x$ to the
    homomorphism $\alpha\mapsto\epsilon(\alpha x)$, where
    $\epsilon:\KzKp\to\Zpl$ is the augmentation, is a homomorphism of
    $A$-modules.  It is an isomorphism since $A$ is the \Zpl-dual of
    the free \Zpl-module \KzKp; see~\cite{CCW05}.  Hence the result
    holds for the case $L=\Zpl$.
    
    For a general \Zpl-module $L$ the natural $A$-module homomorphism 
    $U\Zpl\tensor L\to\Homz(A,L)$ maps into $\Homc(A,L)=UL$.  We can
    define the inverse $A$-homomorphism as follows.  Let 
    $\rho_{j}\in U\Zpl$ be given by $\rho_{j}(\Phi_{k})=\delta_{jk}$.
    If $\sigma\in\Homc(A,L)$ and $\sigma(A_{n})=0$, we send $\sigma$
    to 
    $\sum_{j=0}^{n-1}\rho_{j}\tensor\sigma(\Phi_{j})\in U\Zpl\tensor L$.
\end{proof}

\begin{thm}\label{Uexactness}
    The functor $U$ is exact, and preserves direct sums and direct
    limits.
\end{thm}

\begin{proof}
    As a right adjoint, $U$ is left exact.  On the other hand, if
    $f:L_{1}\to L_{2}$ is a \Zpl-module epimorphism, then
    $Uf:UL_{1}\to UL_{2}$ is an epimorphism of discrete $A$-modules,
    for if $g\in\Homc(A,L_{2})$, then we can define
    $h\in\Homc(A,L_{1})$ with $g=Uf(h)$ as follows.
    For each $n\ge 0$, we can find some $y_{n}\in L_{1}$ such that
    $f(y_{n})=g(\Phi_{n})$, and, since $g(\Phi_{n})=0$ for $n\gg 0$,
    we may choose $y_{n}$ to be~$0$ for $n\gg 0$.  Thus we can define
    $h:A\to L_{1}$ by 
    $h(\sum_{n\ge 0}a_{n}\Phi_{n})=\sum_{n\ge 0}a_{n}y_{n}$.

    That $U$ preserves direct sums is immediate from the definition,
    and an exact functor which preserves direct sums automatically
    preserves direct limits.
\end{proof}

\begin{cor}
    If $D$ is an injective $\Zpl$-module, then $UD$ is injective.
    Hence \DA\ has enough injectives.
\end{cor}

\begin{proof}
    The injectivity of $UD$ follows from adjointness.  As each
    $\Zpl$-module can be embedded in an injective $\Zpl$-module $D$,
    so any discrete $A$-module can be embedded as a $\Zpl$-module in
    such a~$D$.  Then, using adjointness and left-exactness of~$U$,
    any discrete $A$-module can be embedded as an $A$-module in
    some~$UD$.
\end{proof} 

Bousfield introduced in \cite[(7.4)]{Bousfield85} a four-term exact
sequence which underlies the fact, due to Adams and
Baird~\cite{Adams(Baird)}, that all $\Ext^{>2}$ groups are zero.
Bousfield's construction applies to the version of his category
which corresponds to the split summand of \Ktheory; we briefly discuss
this category in the Appendix.  We end this section by constructing the
corresponding exact sequence in~\DA, using the functor~$U$.

\begin{thm}\label{exactsequence}
    For any $M$ in \DA, there is an exact sequence in~\DA:
    $$
    0\to
    M\xrightarrow{\;\alpha\;}
    UM\xrightarrow{\;\beta\;}
    UM\xrightarrow{\;\gamma\;}
    \MQ\to 
    0,
    $$
    where $UM$ denotes the discrete $A$-module obtained by applying
    $U$ to the $\Zpl$-module underlying $M$ (i.e., applying a
    forgetful functor to $M$ before applying~$U$).
\end{thm}

\begin{proof}
    The map $\alpha$ is adjoint to the identity map $M\to M$;
    explicitly $\alpha x$ maps $\theta\in A$ to $\theta x$.  It is
    clear that $\alpha$ is a monomorphism.  For any $x\in M$, the map
    $\alpha x:A\to M$ is $A$-linear.  Moreover, any $A$-linear map
    $f:A\to M$ is determined by~$f(1)$, and
    $f=\alpha\bigl(f(1)\bigr)$.  Thus the image of~$\alpha$ is exactly
    the subset $\Hom_{A}(A,M)$ of $A$-linear maps in $UM=\Homc(A,M)$
    (an $A$-linear map into a discrete $A$-module is automatically
    continuous).

    The map $\beta$ is given by
    $$
    \beta f
    = 
    \Psiq\circ f-f\circ\Psiq,
    $$
    where $f:A\to M$.  Here $\Psiq\circ f$ uses the $A$-module
    structure of~$M$, not that of~$UM$, in other words,
    $\beta f:\theta\mapsto\Psiq f(\theta)-f(\Psiq\theta)$.  It
    is straightforward to check that $\beta$ is an $A$-homomorphism.
    A continuous $\Zpl$-homomorphism from~$A$ into a discrete
    $A$-module which commutes with \Psiq\ is an $A$-module
    homomorphism, hence $\Ker\beta=\Hom_{A}(A,M)=\Im\alpha$.

    \medskip
    
    We define~$\gamma$ as follows.  Let
    $$
    \Theta_{n}^{(j)}(X)
    =
    \prod_{\substack{i=1\\i\ne j}}^{n}(X-q_{i}),
    $$
    and let $\Phi_{n}^{(j)}=\Theta_{n}^{(j)}(\Psiq)\in A$.  Note that
    $\Phi_{n}^{(j)}=\Phi_{n}$ if $j>n$, since
    $\Theta_{n}^{(j)}(X)=\Theta_{n}(X)$ in that case.  If $x\in M$ and
    $\Phi_{n}x=0$, then each $\Phi_{n}^{(j)}x$ is an eigenvector of
    \Psiq\ with eigenvalue~$q_{j}$.  Moreover we have 
    \begin{equation}\label{operatorexp}
        1
        =
        \Phi_{0}
        =
        \sum_{j=1}^{n}\frac{\Phi_{n}^{(j)}}{\Theta_{n}^{(j)}(q_{j})}
    \end{equation}
    in \AQ\ for all~$n\ge 1$.  
    
    If $f\in\Homc(A,M)$, choose $n$ such that $f(\Phi_{k})=0$ for all
    $k\ge n$ and $\Phi_{n}f(\Phi_{k})=0$ for $0\le k<n$.  Then let
    $$
    \gamma f
    =
    \sum_{k\ge 0}
    \sum_{j=1}^{k+1}
    \frac{\Phi_{n}^{(j)}f(\Phi_{k})}
    {\Theta_{n}^{(j)}(q_{j})\Theta_{k+1}^{(j)}(q_{j})}
    \in\MQ.
    $$
    
    Note that we may reverse the order of summation to obtain
    \begin{equation}\label{eigenvectorexp1}
        \gamma f
        =
        \sum_{j=1}^{n}
        \frac{\Phi_{n}^{(j)}x_{j}}{\Theta_{n}^{(j)}(q_{j})},
    \end{equation}
    where
    $$
    x_{j}
    =
    \sum_{k=j-1}^{n-1}
    \frac{f(\Phi_{k})}{\Theta_{k+1}^{(j)}(q_{j})}.
    $$
    In this form it is apparent that the formula for~$\gamma f$ is
    independent of the choice of~$n$.  For if $n\le m$ and 
    $1\le j\le m$, it follows from
    $\Psiq\Phi_{n}^{(j)}x_{j}=q_{j}\Phi_{n}^{(j)}x_{j}$ that
    $$
    \frac{\Phi_{m}^{(j)}x_{j}}{\Theta_{m}^{(j)}(q_{j})}
    =
    \frac{\Phi_{n}^{(j)}x_{j}}{\Theta_{n}^{(j)}(q_{j})}.
    $$
    
    \medskip
    
    Suppose now that $f=\beta g$, with $g(\Phi_{k})=0$ for all $k\ge 
    n$ and $\Phi_{n}g(\Phi_{k})=0$ for $0\le k<n$.  Then
    $$
    x_{j}
    =
    \sum_{k=j-1}^{n-1}
    \frac{(\Psiq-q_{k+1})g(\Phi_{k})-g(\Phi_{k+1})}
    {\Theta_{k+1}^{(j)}(q_{j})}.
    $$
    But, since
    $\Theta_{k+1}^{(j)}(q_{j})=(q_{j}-q_{k+1})\Theta_{k}^{(j)}(q_{j})$
    for $1\le j\le k$, and $g(\Phi_{n})=0$,
    $$    
    x_{j}
    =
    \sum_{k=j-1}^{n-1}
    \frac{(\Psi^{q}-q_{j})g(\Phi_{k})}{\Theta_{k+1}^{(j)}(q_{j})}.
    $$
    It follows from~\eqref{eigenvectorexp1} that $\gamma\beta g=0$,
    since $\Phi_{n}^{(j)}(\Psi^{q}-q_{j})=\Phi_{n}$.
    
    \medskip
    
    To prove that $\Ker\gamma\subseteq\Im\beta$ we need first the 
    following lemma.
    
    \begin{lemma}\label{torsion-hom}
        Suppose $f\in\Ker\gamma$, where $M$ is a discrete $A$-module
        and $n$ is chosen as above.  Then
        $\sum_{k=0}^{n-1}(\Phi_n/\Phi_{k+1})f(\Phi_{k})$ is a
        \Zpl-torsion element of~$M$.
    \end{lemma}

    \begin{proof}
        Since $\gamma f=0$, in equation~\eqref{eigenvectorexp1} each
        $\Phi_{n}^{(j)}x_{j}=0$ in \MQ\ , since otherwise \Psiq\ would
        have linearly dependent eigenvectors with distinct eigenvalues.
        Thus for each $j=1,\dots,n$,
        \begin{equation}\label{keysum-hom}
            0
            =
            \sum_{k=j-1}^{n-1}
            \frac{\Phi_{n}^{(j)}f(\Phi_{k})}
            {\Theta_{k+1}^{(j)}(q_{j})}.
        \end{equation}
    
        Now if $j\le k+1\le n$,
        $\Phi_{n}^{(j)}=\Phi_{k+1}^{(j)}(\Phi_{n}/\Phi_{k+1})$ in~$A$.
        Hence \eqref{keysum-hom} becomes
        $$
        0
        =
        \sum_{k=j-1}^{n-1}
        \frac{\Phi_{k+1}^{(j)}(\Phi_{n}/\Phi_{k+1})f(\Phi_{k})}
        {\Theta_{k+1}^{(j)}(q_{j})}
        $$
        in \MQ. Let $d_{n}$ denote the least common multiple of the
        $\Theta_{k+1}^{(j)}(q_{j})$ for $j\le k+1\le n$.  Then
        multiplying by $d_{n}$ yields 
        $0=z_{j}$ in \MQ, where
        $$
        z_{j}
        :=
        \sum_{k=j-1}^{n-1}
        \bigl(
	\bigl(d_{n}/\Theta_{k+1}^{(j)}(q_{j})\bigr)\Phi_{k+1}^{(j)}
        \bigr)
        (\Phi_{n}/\Phi_{k+1})f(\Phi_{k})
        \in M.
        $$
        Hence $z_{j}$ is a \Zpl-torsion element in~$M$, and so is
        \begin{align*}
	    \sum_{j=1}^{n}z_{j}
	    &=
	    \sum_{k=0}^{n-1}\left(
	    \sum_{j=1}^{k+1}
	    \bigl(d_{n}/\Theta_{k+1}^{(j)}(q_{j})\bigr)\Phi_{k+1}^{(j)}
	    \right)
	    (\Phi_{n}/\Phi_{k+1})f(\Phi_{k})\\
	    &=
	    d_{n}\sum_{k=0}^{n-1}
	    (\Phi_{n}/\Phi_{k+1})f(\Phi_{k}),
        \end{align*}
        where we use~\eqref{operatorexp}. The result follows.
    \end{proof}

    Now if $f\in\Homc(A,M)$ and $r\ge 1$, define
    $\ft_{r}\in\Homc(A,M)$ by
    $$
    \ft_{r}(\Phi_{k})
    =
    \begin{cases}
        \displaystyle
        \sum_{i=0}^{k-1}(\Phi_{k}/\Phi_{i+1})f(\Phi_{i}),& 
        \text{if $1\le k\le r$,}\\[2.4ex]
        0,&
        \text{if $k=0$ or $k>r$.}
    \end{cases}
    $$
    It is a simple calculation that
    $$
    (\beta\ft_{r}+f)(\Phi_{k})
    =
    \begin{cases}
        0,&
        \text{if $k<r$,}\\[0.6ex]
        \displaystyle\sum_{i=0}^{r}(\Phi_{r+1}/\Phi_{i+1})f(\Phi_{i}),& 
        \text{if $k=r$,}\\[2.4ex]
        f(\Phi_{k}),&
        \text{if $k>r$.}
    \end{cases}
    $$
    Suppose now that $f\in\Ker\gamma$, that $f(\Phi_{k})=0$ for all
    $k\ge n$ and $\Phi_{n}f(\Phi_{k})=0$ for $0\le k<n$.  Then if
    $r>n$,
    $$
    (\beta\ft_{r}+f)(\Phi_{r})
    =
    (\Phi_{r+1}/\Phi_{n})y,
    $$
    where $y=\sum_{i=0}^{n-1}(\Phi_{n}/\Phi_{i+1})f(\Phi_{i})$ is, by
    Lemma~\ref{torsion-hom}, a \Zpl-torsion element of~$M$.  But now
    Lemma~\ref{PhiQuotientsModp} shows that for $r$ sufficiently large
    $(\Phi_{r+1}/\Phi_{n})y=0$, in which case $f=\beta(-\ft_{r})$.
    Hence we have shown that $\Ker\gamma\subseteq\Im\beta$.
    
    \medskip
    
    It remains to show that $\gamma$ maps onto \MQ. Let $x\in M$, and
    suppose $\Phi_{n}x=0$.  Letting $d$ denote the least common
    multiple of the $\Theta_{n}^{(k)}(q_{k})$ for $k=1,\dots,n$,
    define $f\in\Homc(A,M)$ by 
    $f(\Phi_{k}) 
    =
    \bigl(d\Theta_{k+1}^{(k)}(q_{k})/\Theta_{n}^{(k)}(q_{k})\bigr)
    \Phi_{n}^{(k)}x$.
    Since
    $\Phi_{n}^{(j)}\Phi_{n}^{(k)}x=\Theta_{n}^{(j)}(q_{k})\Phi_{n}^{(k)}x$
    is zero unless $j=k$, it is a simple calculation
    using~\eqref{operatorexp} that $\gamma f=dx$ in \MQ.
    
    The proof that $\gamma$ is an epimorphism is completed by showing
    that if $f\in\Homc(A,M)$, there exists $g\in\Homc(A,M)$ such that
    $f-pg\in\Ker\gamma$.

    To see this, choose a multiple~$m$ of~$2p-2$ such that
    $f(\Phi_{k})=0$ for all $k\ge m$ and $\Phi_{m}f(\Phi_{k})=0$ for
    all $k\ge 0$.  By Lemma~\ref{PhiQuotientsModp}, for each $k\ge m$
    there is a $\theta_{k}\in A$ such that
    $(\Phi_{k+1}/\Phi_{k-m+1})=\Phi_{m}+p\theta_{k}$.  We define $g$
    as follows
    $$
    g(\Phi_{k})
    =
    \begin{cases}
        \theta_{k}f(\Phi_{k-m}),&
        \text{if $k\ge m$,}\\
        0,& 
        \text{otherwise.}
    \end{cases}
    $$

    For suitably large~$n$, let
    $$
    y_{j}
    =
    \sum_{k=j-1}^{n-1}
    \frac{pg(\Phi_{k})}{\Theta_{k+1}^{(j)}(q_{j})}
    =
    \sum_{k=j-m-1}^{n-m-1}
    \frac{(\Phi_{k+m+1}/\Phi_{k+1})f(\Phi_{k})}{\Theta_{k+m+1}^{(j)}(q_{j})},
    $$
    so that $p\gamma
    g=\sum_{j=1}^{n}\Phi_{n}^{(j)}y_{j}/\Theta_{n}^{(j)}(q_{j})$.  It
    is now clear that $p\gamma g=\gamma f$.  For
    $$
    (\Phi_{k+m+1}/\Phi_{k+1})\Phi_{n}^{(j)}f(\Phi_{k})
    =
    \left(
    \prod_{i=k+2}^{k+m+1}\!(q_{j}-q_{i})
    \right)
    \Phi_{n}^{(j)}f(\Phi_{k}),
    $$
    which is zero if $k+2\le j\le k+m+1$, while
    $$
    \Theta_{k+m+1}^{(j)}(q_{j})
    =
    \Theta_{k+1}^{(j)}(q_{j})
    \!\prod_{i=k+2}^{k+m+1}\!(q_{j}-q_{i}),
    $$
    if $1\le j\le k+1$.
\end{proof}

\section{Bousfield's Main Result}
\label{Bresult}

We translate Bousfield's main result, giving an algebraic
classification of $K_{(p)}$-local homotopy types, into our language of
discrete $A$-modules.  We claim no originality here, but we believe it
is useful to summarise Bousfield's results in our language.

There are several steps in Bousfield's construction: we need a graded
version of our discrete module category, we need to understand $\Ext$
groups in this category, and we need to see $k$-invariants associated
to $K_{(p)}$-local spectra as elements in such $\Ext$ groups.  We
outline these steps without proof.  The proofs can be easily adapted
from those in~\cite{Bousfield85}.

\newcommand{\myitem}[1]{\item[#1]\quad\par\smallskip}
\newcounter{n}
\begin{list}{}%
{\usecounter{n}%
 \setlength{\topskip}{\bigskipamount}
 \setlength{\itemsep}{\medskipamount}
 \setlength{\labelwidth}{0pt}
 \setlength{\leftmargin}{0pt}
 \renewcommand{\makelabel}[1]{\stepcounter{n}\arabic{n}) #1.}}

\myitem{The category $\DA_{*}$}
For $i\in\Z$, there is an automorphism $T^{i}:\DA\to\DA$ with $T^{i}M$
equal to $M$ as a $\Zpl$-module but with $\Psiq:T^{i}M\to T^{i}M$ equal to
$q^{i}\Psiq:M\to M$.  An object of $\DA_{*}$ is a collection of objects
$M_{n}\in\DA$ for $n\in\Z$ together with isomorphisms 
$u:TM_{n}\cong M_{n+2}$ in~\DA\ for all~$n$.  A morphism $f:M\to N$ in
$\DA_{*}$ is a collection of morphisms $f_{n}:M_{n}\to N_{n}$ in~\DA\
such that $uf_{n}=f_{n+2}u$ for all $n\in\Z$.

The point of this construction is that $K_{*}(X;\Zpl)$ is an object of
$\DA_{*}$ for any spectrum $X$.

\myitem{$\Ext$ groups in $\DA_{*}$}
The category $\DA_{*}$ has enough injectives, allowing the definition of
(bigraded) $\Ext$ groups.  The groups $\Ext^{s,t}_{\DA_{*}}(\_,\_)$ vanish
for $s>2$, essentially as a consequence of the exact sequence of
Theorem~\ref{exactsequence}.  There is an Adams spectral sequence
with
$$
E_{2}^{s,t}(X,Y)
=
\Ext^{s,t}_{\DA_{*}}\bigl(K_{(p)*}(X),K_{(p)*}(Y)\bigr),
$$
converging strongly to $[X_{K_{(p)}}, Y_{K_{(p)}}]_{*}$.

\myitem{$k$-invariants and the category $k\DA_{*}$} To each
$K_{(p)}$-local spectrum $X$ is associated a $k$-invariant
$k_{X}\in\Ext_{\DA_{*}}^{2,1}\bigl(K_{{(p)}*}(X),
K_{{(p)}*}(X)\bigr)$.  The only non-trivial differential $d_{2}$ in
the Adams spectral sequence can be expressed in terms of these
$k$-invariants.  We form the additive category $k\DA_{*}$ whose
objects are pairs $(M,\kappa)$, with $M\in\DA_{*}$ and
$\kappa\in\Ext_{\DA_{*}}^{2,1}(M,M)$, and whose morphisms from
$(M,\kappa)$ to $(N,\lambda)$ are morphisms $f:M\to N$ in $\DA_{*}$
with $\lambda f=f\kappa\in\Ext_{\DA_{*}}^{2,1}\!(M,N)$.

\end{list}

This now allows us to translate the main result of~\cite{Bousfield85}
into our setting.

\begin{thm}\label{localobjects}
    Homotopy types of $K_{(p)}$-local spectra are in one-to-one
    correspondence with isomorphism classes in $k\DA_{*}$.  
    \qed
\end{thm}

\section{Appendix: The Split Setting}
    
In this section we summarise the Adams summand analogues of our
results.  We omit proofs as these are easy adaptations of those given
in the preceding sections.

For a fixed odd prime $p$, we denote the periodic Adams summand
by~$G$, and we write $B$ for the ring of stable degree zero operations
$G^{0}(G)$.  As usual, we choose $q$ primitive modulo $p^{2}$, and we
set $\qh=q^{p-1}$.  Let
$\Phih_{n}=\prod_{i=1}^{n}(\Psiq-\qh^{(-1)^{i}\floor{i/2}})\in B$.
The elements $\Phih_{n}$, for $n\ge 0$, form a topological \Zpl-basis
for~$B$; see~\cite[Theorem 6.13]{CCW05}.  Just as for~$A$, the
topological ring $B$ can be viewed as a completion of the polynomial
ring $\Zpl[\Psiq]$.  The results in the two cases are formally very
similar, differing only in that in many formulas $q$ must be replaced
by~$\qh$.

Let \DB\ denote the category of discrete $B$-modules, defined in the
obvious way by analogy with Definition~\ref{discrete}.  There is a
corresponding graded version $\DB_{*}$.  Then the additive category
$k\DB_{*}$ is formed just as we did in the non-split setting.  Its
objects are pairs $(M,\kappa)$, where $M\in\DB_{*}$ and
$\kappa\in\Ext_{\DB_{*}}^{2,1}\!(M,M)$.

\begin{thm}\quad\par
    \begin{enumerate}
        \item \DB\ is a cocomplete abelian category with enough
        injectives.

        \item For any spectrum $X$, $G_{0}(X)$ is an object of \DB,
        and $G_{*}(X)$ is an object of $\DB_{*}$.
    
        \item \DB\ is isomorphic to Bousfield's category
        $\mathcal{B}(p)$.
    
        \item The functor $U(\_)=\Homc(B,\_)$ from \Zpl-modules to \DB\
        is right adjoint to the forgetful functor.
    
        \item For any $M$ in \DB, there is an exact sequence in
        \DB
        $$
        0\to
        M\xrightarrow{\;\alpha\;}
        UM\xrightarrow{\;\beta\;}
        UM\xrightarrow{\;\gamma\;}
        \MQ\to 
        0.
        $$
    
        \item The groups $\Ext^{s,t}_{\DB_{*}}\,(\_,\_)$ vanish for
        $s>2$.
    
        \item Homotopy types of $G$-local spectra are in one-to-one
        correspondence with isomorphism classes in $k\DB_{*}$.  
        \qed
    \end{enumerate}
\end{thm}

\end{document}